%% file: paper.tex
\documentclass[11pt]{article}
\usepackage[top=1in, bottom=1in, left=1in, right=1in]{geometry}
\usepackage{tikz}
\usetikzlibrary{shapes,arrows,matrix,patterns,positioning}
\usetikzlibrary{calc}
\usepackage{color}
\usepackage{graphicx}
\usepackage{algorithmic}
\usepackage{amssymb}
\usepackage{epstopdf}
\usepackage[ruled]{algorithm2e}
\usepackage{float}
\usepackage{amsmath,amsthm,bm,color,epsfig,enumerate,caption}
\usepackage{hyperref}
\usepackage{xparse}

\newtheorem{theorem}{Theorem}[section]

\newcommand{\norm}[1]{\left\lVert#1\right\rVert}
\renewcommand{\L}{\mathcal{L}}

\newcommand{\dist}{\text{dist}}

\renewcommand{\O}{\mathcal{O}}

\newcommand{\eps}{\epsilon}
\newcommand{\p}{\partial}
\newcommand{\R}{\mathbb{R}}
\renewcommand{\i}{\imath}

\tikzset{%
    Cote node/.style={%
        midway,
        sloped,
        fill=white,
        inner sep=1.5pt,
        outer sep=2pt
    },
    Cote arrow/.style={%
        <->,
        >=latex,
        very thin
    }
}

\makeatletter
\newcommand*{\extendadd}{%
  \mathbin{%
    \mathpalette\extend@add{}%
  }%
}
\newcommand*{\extend@add}[2]{%
  \ooalign{%
    $\m@th#1\leftrightarrow$%
    \vphantom{$\m@th#1\updownarrow$}
    \cr
    \hfil$\m@th#1\updownarrow$\hfil
  }%
}

\NewDocumentCommand{\Cote}{%
    s       
    D<>{1.5pt} 
    O{.75cm}    
    m       
    m       
    m       
    D<>{o}  
    O{}     
    }{%

    {\tikzset{#8}

    \coordinate (@1) at #4 ;
    \coordinate (@2) at #5 ;

    \if #7v 
        \coordinate (@0) at ($($#4!.5!#5$) + (#3,0)$) ;
        \coordinate (@4) at (@0|-@1) ;
        \coordinate (@5) at (@0|-@2) ;
    \else
    \if #7h 
        \coordinate (@0) at ($($#4!.5!#5$) + (0,#3)$) ;
        \coordinate (@4) at (@0-|@1) ;
        \coordinate (@5) at (@0-|@2) ;
    \else 
    \ifnum\pdfstrcmp{\unexpanded\expandafter{\@car#7\@nil}}{(}=\z@
        \coordinate (@5) at ($#7!#3!#5$) ;
        \coordinate (@4) at ($#7!#3!#4$) ;
    \else 
        \coordinate (@5) at ($#5!#3!90:#4$) ;
        \coordinate (@4) at ($#4!#3!-90:#5$) ;
    \fi\fi\fi

    \draw[very thin,shorten >= #2,shorten <= -2*#2] (@4) -- #4 ;
    \draw[very thin,shorten >= #2,shorten <= -2*#2] (@5) -- #5 ;

    \IfBooleanTF #1 {
    \draw[Cote arrow,-] (@4) -- (@5)
        node[Cote node] {#6\strut};
    \draw[Cote arrow,<-] (@4) -- ($(@4)!-6pt!(@5)$) ;
    \draw[Cote arrow,<-] (@5) -- ($(@5)!-6pt!(@4)$) ;
    }{
    \ifnum\pdfstrcmp{\unexpanded\expandafter{\@car#7\@nil}}{(}=\z@
        \draw[Cote arrow] (@5) to[bend right]
            node[Cote node] {#6\strut} (@4) ;
    \else
    \draw[Cote arrow] (@4) -- (@5)
        node[Cote node] {#6\strut};
    \fi
    }}
}
\makeatother

\begin{document}

\title{A Multiscale Butterfly Algorithm for Multidimensional Fourier
  Integral Operators}

\author{Yingzhou Li$^{\sharp}$, Haizhao Yang$^{\dagger}$\footnote{Corresponding author. Email address: haizhao@math.stanford.edu.} and Lexing Ying$^{\dagger\sharp}$\\
  \vspace{0.1in}\\
  $\dagger$ Department of Mathematics, Stanford University\\
  $\sharp$ ICME, Stanford University
}

\date{Nov 2014}
\maketitle

\begin{abstract}
  This paper presents an efficient multiscale butterfly algorithm for
  computing Fourier integral operators (FIOs) of the form $(\L f)(x) =
  \int_{\mathbb{R}^d}a(x,\xi) e^{2\pi \i \Phi(x,\xi)}\widehat{f}(\xi)
  d\xi$, where $\Phi(x,\xi)$ is a phase function, $a(x,\xi)$ is an
  amplitude function, and $f(x)$ is a given input. The frequency
  domain is hierarchically decomposed into a union of Cartesian
  coronas. The integral kernel $a(x,\xi) e^{2\pi \i \Phi(x,\xi)}$ in
  each corona satisfies a special low-rank property that enables the
  application of a butterfly algorithm on the Cartesian phase-space
  grid. This leads to an algorithm with quasi-linear operation
  complexity and linear memory complexity. Different from previous
  butterfly methods for the FIOs, this new approach is simple and
  reduces the computational cost by avoiding extra coordinate
  transformations. Numerical examples in two and three dimensions are
  provided to demonstrate the practical advantages of the new
  algorithm.
\end{abstract}



\section{Introduction}
\label{sec:intro}

This paper is concerned with the rapid application of Fourier
integral operators (FIOs), which are defined as
\begin{equation}
  \label{eqn:FIO}
  (\L f)(x) = \int_{\mathbb{R}^d}a(x,\xi)  e^{2\pi \i \Phi(x,\xi)}\widehat{f}(\xi) d\xi,
\end{equation}
where
\begin{itemize}
\item $a(x,\xi)$ is an amplitude function that is smooth both in
  $x$ and $\xi$,
\item $\Phi(x,\xi)$ is a phase function that is smooth in $(x,\xi)$
  for $\xi\neq 0$ and obeys the homogeneity condition of degree $1$ in
  $\xi$, namely, $\Phi(x,\lambda \xi)=\lambda\Phi(x,\xi)$ for each
  $\lambda>0$;
\item $\widehat{f}$ is the Fourier transform of the input $f$ defined
  by
  \[
  \widehat{f}(\xi) = 
  \int_{\mathbb{R}^d} e^{-2\pi \i x \cdot \xi} f(x) dx.
  \]
\end{itemize}
The computation of Fourier integral operators appears quite often
in the numerical solution of wave equations and related
applications in computational geophysics. In a typical setting, it
is often assumed that the problem is periodic (i.e., $a(x,\xi)$,
$\Phi(x,\xi)$, and $f(x)$ are all periodic in $x$) or the function
$f(x)$ decays sufficiently fast so that one can embed the problem in a
sufficiently large periodic cell. A simple discretization in two
dimensions considers functions $f$ given on a Cartesian grid
\begin{equation}
  \label{eqn:X}
  X = \left\{ x = \left( \frac{n_1}{N}, \frac{n_2}{N}\right), 0 \leq  n_1,n_2 < N\text{ with }
  n_1, n_2 \in \mathbb{Z} \right\}
\end{equation}
in a unit square and defines the discrete Fourier integral operator
by
\[
(Lf)(x) = \sum_{\xi\in\Omega} a(x,\xi) e^{2\pi \i \Phi(x,\xi)}\widehat{f}(\xi), \quad x\in X,
\]
where
\begin{equation}
  \label{eqn:Omega}
  \Omega = \left\{ \xi = (n_1, n_2),- \frac{N}{2} \leq n_1,n_2 < \frac{N}{2}\text{ with }
  n_1, n_2 \in \mathbb{Z} \right\},
\end{equation}
and $\widehat{f}$ is the discrete Fourier transform of $f$
\[
\widehat{f}(\xi) = \frac{1}{N^2} \sum_{x\in X} e^{-2\pi \i x \cdot \xi} f(x).
\]

In most examples, since $a(x,\xi)$ is a smooth symbol of order zero
and type $(1,0)$ \cite{symbol1,FIO07,symbol2,Theory,symbol3},
$a(x,\xi)$ is numerically low-rank in the joint $X$ and $\Omega$
domain and its numerical treatment is relatively easy. Therefore, we
will simplify the problem by assuming $a(x,\xi)=1$ in the following
analysis and the algorithm description. \cite{FIO09} is referred to for
discussion on how to deal with a non-constant amplitude
function. Under this assumption, the discrete FIO discussed in this
paper takes the following form:
\begin{equation}
  \label{eqn:FIOdiscrete}
  (L f)(x)= \sum_{\xi \in \Omega} e^{2\pi \i \Phi(x,\xi)} \widehat{f}(\xi), \quad x\in X.
\end{equation}
A direct computation of \eqref{eqn:FIOdiscrete} takes $\O(N^4)$
operations, which is quadratic in the number of DOFs $N^2$.  Hence, a
practical need is to design efficient and accurate algorithms to
evaluate \eqref{eqn:FIOdiscrete}. This research topic is of great
interest for computing wave equations especially in geophysics
\cite{Yingwave,Hu,invRadon,invFIO}.

\subsection{Previous work}
An earlier method for the rapid computation of general FIOs is the
algorithm for two dimensional problems proposed in \cite{FIO07}. This
method starts by partitioning the frequency domain $\Omega$ into
$\O(\sqrt{N})$ wedges of equal angle. The integral
\eqref{eqn:FIOdiscrete} {\em restricted to each wedge} is then
factorized into two components, both of which can be handled efficiently.
The first one has a low-rank structure that leads to an $\O(N^2\log
N)$ fast computation, while the second one is a non-uniform Fourier
transform which can be evaluated in $\O(N^2\log N)$ steps with the
algorithms developed in \cite{NFFT1,NFFT2,NFFT3}. Summing the
computational cost over all $\O(\sqrt{N})$ wedges gives an
$\O(N^{2.5}\log N)$ computational cost.

Shortly after, an algorithm with quasilinear complexity for general FIOs was
proposed in \cite{FIO09} using the framework of the butterfly
algorithms in \cite{Butterfly1,Butterfly2}. This approach introduces a
polar coordinate transformation in the frequency domain to remove the
singularity of $\Phi(x,\xi)$ at $\xi=0$, proves the existence of
low-rank separated approximations between certain pairs of spatial and
frequency domains, and implements the low-rank approximations with
oscillatory Chebyshev interpolations. The resulting algorithm
evaluates \eqref{eqn:FIOdiscrete} with $\O(N^2\log N)$ operations and
$\O(N^2)$ memory, both essentially linear in terms of the number of
unknowns.

Another related research direction seeks for sparse representations of
the FIOs using modern basis functions from harmonic analysis. A sparse
representation allows fast matrix-vector products in the transformed
domain. Local Fourier transforms \cite{Basis1,Basis2,Basis5},
wavelet-packet transforms \cite{Basis4}, the curvelet transform
\cite{Candes2006,Candes2004,CandesI,CandesII}, and the wave atom frame
\cite{Demanet2007,Basis3} have been investigated for the purpose of
operator sparsification. In spite of favorable asymptotic behaviors,
the actual representations of the FIOs typically have a large
pre-factor constant in terms of both the computational time and the
memory requirement. This makes them less competitive compared to the
approaches in \cite{FIO07,FIO09}.

\subsection{Motivation}

The main motivation of the current work is to improve the performance
of the butterfly algorithm in \cite{FIO09}. As we pointed out earlier,
this algorithm starts by applying a polar coordinate transformation in
the frequency domain to remove the singularity of the phase function
at $\xi=0$. For this reason, we refer to this algorithm as the {\bf
  polar Butterfly algorithm}. More precisely, the polar butterfly
algorithm introduces a polar-Cartesian coordinate transformation
$T:(p_1,p_2)\rightarrow (\xi_1,\xi_2)$ such that
\begin{equation}
  \label{eqn:TF}
  \xi=(\xi_1,\xi_2) = \frac{\sqrt{2}}{2}Np_1e^{2\pi \i p_2},
  \quad e^{2\pi \i p_2}=(\cos 2\pi p_2,\sin 2\pi p_2).
\end{equation}
Let $P=T^{-1}(\Omega)$. By definition, each point $p=(p_1,p_2)\in
P$ belongs to $[0,1]^2$. The new phase function $\Psi(x,p)$ in the $p$
variable is now given by
\begin{equation}
  \label{eqn:newPhase}
  \Psi(x,p):=\frac{1}{N}\Phi(x,\xi)
  = \frac{\sqrt{2}}{2}\Phi(x,e^{2\pi\i p_2})p_1,
\end{equation}
where the last identity comes from the homogeneity of $\Phi(x,\xi)$ in
$\xi$.  Thus, computing \eqref{eqn:FIOdiscrete} is equivalent to
evaluate
\begin{equation}
  \label{eqn:09}
  (L f)(x)= \sum_{\xi \in \Omega}  e^{2\pi \i \Phi(x,\xi)}\widehat{f}(\xi)
  = \sum_{p \in P} e^{2\pi \i N\Psi(x,p)}\widehat{f}(T(p)).
\end{equation}
The new phase function $\Psi(x,p)$ is smooth in the whole domain
$(x,p)\subset [0,1]^2\times[0,1]^2$, since $\Phi(x,\xi)$ is smooth in
$(x,\xi)$ for $\xi\neq 0$. This smoothness guarantees a low-rank
separated approximation of $e^{2\pi \i N\Psi(x,p)}$ when $x$ and $p$
are properly restricted to certain subdomains in $X\times P$ under
certain geometric configuration. This low-rank property allows for the
application of the butterfly algorithm in \cite{Butterfly5} and
results in a fast algorithm with an $\O(N^2\log N)$ computational
complexity and an $\O(N^2)$ memory complexity.

However, the application of this polar-Cartesian transformation comes
with several drawbacks, which result in a large pre-factor of the
computational complexity. First, due to the polar grid in the
frequency domain, the points in $P$ for the butterfly algorithm are
irregularly distributed and a separate Chebyshev interpolation matrix
is required for the evaluation at each point. In order to avoid the
memory bottleneck from storing these interpolation matrices, the polar
butterfly algorithm generates these interpolation matrices on-the-fly
during the evaluation. This turns out to be expensive in the operation
count. Second, since the amplitude and phase functions are often
written in the Cartesian coordinates, the polar butterfly algorithm
applies the polar-Cartesian transformation for each kernel
evaluation. Finally, in order to maintain a reasonable accuracy, the
polar butterfly algorithm divides the frequency domain into multiple
parts and applies the same butterfly algorithm to each part
separately. This also increases the actual running time by a
non-trivial constant factor.

\subsection{Our contribution}

Those drawbacks of the polar Butterfly algorithm motivate us to
propose a multiscale butterfly algorithm using a Cartesian grid both
in the spatial and frequency domain.  To deal with the singularity of
the kernel $\Phi(x,\xi)$ at $\xi=0$, we hierarchically decompose the
frequency domain into a union of non-overlapping Cartesian coronas
with a common center $\xi=0$ (see Figure \ref{fig:domain-decomp}). More
precisely, define
\[
\Omega_j = \left\{(n_1,n_2):  \frac{N}{2^{j+1}} < \max(|n_1|,|n_2|) \le \frac{N}{2^j} \right\} \cap \Omega
\]
for $j=1,\dots,\log N-s$, where $s$ is just a small constant integer.
The domain $\Omega_d = \Omega\setminus \cup_j\Omega_j$ is the
remaining square grid at the center of constant size. Following this
decomposition of the frequency domain, one can write
\eqref{eqn:FIOdiscrete} accordingly as
\begin{equation}
  \label{eqn:HFIO}
  (L f)(x)=
  \sum_j \left( \sum_{\xi \in \Omega_j} e^{2\pi \i \Phi(x,\xi)}\widehat{f}(\xi) \right) +
  \sum_{\xi \in \Omega_d} e^{2\pi \i \Phi(x,\xi)}\widehat{f}(\xi).
\end{equation}

\input{fig-domain-decomp}

The kernel function of \eqref{eqn:HFIO} is smooth in each sub-domain
$\Omega_j$ and a Cartesian butterfly algorithm is applied to evaluate
the contribution from $\Omega_j$.  For the center square $\Omega_d$,
since it contains only a constant number of points, a direct summation
is used. Because of the mutiscale nature of the frequency domain
decompositions, we refer to this algorithm as {\bf the multiscale
  butterfly algorithm}.  As we shall see, the computational and memory
complexity of the multiscale butterfly algorithm are still $\O(N^2\log
N)$ and $\O(N^2)$, respectively. On the other hand, the pre-factors
are much smaller, since the multiscale butterfly is based on the
Cartesian grids and requires no polar-Cartesian transformation.

\subsection{Organization}

The rest of this paper is organized as follows. Section
\ref{sec:butterfly} presents the overall structure of a butterfly
algorithm. Section \ref{sec:LRapproximations} proves a low-rank
property that is essential to the multiscale butterfly
algorithm. Section \ref{sec:HBA} combines the results of the previous
two sections and describes the multiscale butterfly algorithm in
detail. In Section \ref{sec:results}, numerical results of several
examples are provided to demonstrate the efficiency of the multiscale
butterfly algorithm. Finally, we conclude this paper with some
discussion in Section \ref{sec:conclusion}.

\section{The Butterfly Algorithm}
\label{sec:butterfly}

This section provides brief description of the overall structure of
the butterfly algorithm. In this section, $X$ and $\Omega$ refer to
two general sets of $M$ points in $\R^2$, respectively. We assume the
points in these two sets are distributed quasi-uniformly but they are
not necessarily the sets defined in \eqref{eqn:X} and
\eqref{eqn:Omega}.

Given an input $\{g(\xi), \xi \in \Omega\}$, the goal is to compute
the potentials $\{u(x), x\in X\}$ defined by
\[
u(x) = \sum_{\xi\in \Omega} K(x,\xi) g(\xi), \quad x \in X,
\]
where $K(x,\xi)$ is a kernel function. Let $D_X \supset X$ and
$D_\Omega \supset \Omega$ be two square domains containing $X$ and
$\Omega$ respectively. The main data structure of the butterfly
algorithm is a pair of quadtrees $T_X$ and $T_\Omega$. Having $D_X$ as
its root box, the tree $T_X$ is built by recursive dyadic partitioning
of $D_X$ until each leaf box contains only a few points. The tree
$T_\Omega$ is constructed by recursively partitioning in the same
way. With the convention that a root node is at level 0, a leaf node
is at level $L = O(\log M)$ under the quasi-uniformity condition about
the point distributions, where $M$ is the number of points in $X$ and
$\Omega$.  Throughout, we shall use $A$ and $B$ to denote the square
boxes of $T_X$ and $T_\Omega$ with $\ell_A$ and $\ell_B$ denoting
their levels, respectively.

At the heart of the butterfly algorithm is a special low-rank
property. Consider any pair of boxes $A \in T_X$ and $B \in T_\Omega$
obeying the condition $\ell_A + \ell_B = L$. The butterfly algorithm
assumes that the submatrix $\{K(x,\xi)\}_{x \in A, \xi \in B}$ to be
approximately of a constant rank. More precisely, for any $\eps$,
there exists a constant $r_\eps$ independent of $M$ and two sets of
functions $\{\alpha^{AB}_t (x)\}_{1\le t \le r_\eps}$ and
$\{\beta^{AB}_t (\xi)\}_{1\le t \le r_\eps}$ such that the following
holds
\begin{equation}
  \left| K(x,\xi) - \sum_{t=1}^{r_\eps} \alpha^{AB}_t(x)
  \beta^{AB}_t(\xi) \right| \le \eps, \quad \forall x\in A, \forall
  \xi\in B.
  \label{eqn:glr}
\end{equation}
The number $r_\eps$ is called the {\em $\eps$-separation rank}. The
exact form of the functions $\{\alpha^{AB}_t (x)\}_{1\le t\le r_\eps}$
and $\{\beta^{AB}_t (\xi)\}_{1\le t \le r_\eps}$ of course depends on
the problem to which the butterfly algorithm is applied.

For a given square $B$ in $D_\Omega$, define $u^B(x)$ to be the {\em
  restricted potential} over the sources $\xi\in B$
\[
u^B(x) = \sum_{\xi\in B} K(x,\xi) g(\xi).
\]
The low-rank property gives a compact expansion for $\{u^B(x)\}_{x\in
  A}$ as summing \eqref{eqn:glr} over $\xi\in B$ with weights $g(\xi)$
gives
\[
\left| u^B(x) - \sum_{t=1}^{r_\eps} \alpha^{AB}_t(x) \left( \sum_{\xi\in B} \beta^{AB}_t(\xi) g(\xi) \right) \right|
\le \left( \sum_{\xi\in B} |g(\xi)| \right) \eps,
\quad \forall x \in A.
\]
Therefore, if one can find coefficients $\{\delta^{AB}_t\}_{1\le t\le
  r_\eps}$ obeying
\begin{equation}
  \delta^{AB}_t \approx \sum_{\xi\in B} \beta^{AB}_t(\xi) g(\xi), \quad 1\le t\le r_\eps,
  \label{eqn:delta}
\end{equation}
then the restricted potential $\{u^B(x)\}_{x\in A}$ admits a compact
expansion
\[
\left| u^B(x) - \sum_{t=1}^{r_\eps} \alpha^{AB}_t(x) \delta^{AB}_t \right| \le \left( \sum_{\xi\in B} |g(\xi)| \right) \eps,
\quad \forall x\in A.
\]
A key point of the butterfly algorithm is that for each pair $(A,B)$,
the number of terms in the expansion is independent of $M$.

Computing $\{\delta^{AB}_t\}_{1\le t \le r_\eps}$ by means of
\eqref{eqn:delta} for all pairs $A, B$ is not efficient when $B$ is a
large box because for each $B$ there are many paired boxes $A$.  The
butterfly algorithm, however, comes with an efficient way for
computing $\{\delta^{AB}_t\}_{1\le t\le r_\eps}$ recursively.  The
general structure of the algorithm consists of a top down traversal of
$T_X$ and a bottom up traversal of $T_\Omega$, carried out
simultaneously.
\begin{enumerate}
\item Construct the trees $T_X$ and $T_\Omega$ with root nodes $D_X$
  and $D_\Omega$.

\item Let $A$ be the root of $T_X$. For each leaf box $B$ of
  $T_\Omega$, construct the expansion coefficients $\{
  \delta^{AB}_t\}_{1\le t \le r_\eps}$ for the potential
  $\{u^B(x)\}_{x\in A}$ by simply setting
  \begin{equation}
    \delta^{AB}_t = \sum_{\xi\in B} \beta^{AB}_t(\xi) g(\xi), \quad 1\le t \le r_\eps.
  \label{eqn:bf1}
  \end{equation}

\item For $\ell = 1, 2, \ldots, L$, visit level $\ell$ in $T_X$ and
  level $L-\ell$ in $T_\Omega$. For each pair $(A,B)$ with $\ell_A =
  \ell$ and $\ell_B = L-\ell$, construct the expansion coefficients
  $\{\delta^{AB}_t\}_{1\le t \le r_\eps}$ for the potential
  $\{u^B(x)\}_{x\in A}$ using the low-rank representation constructed
  at the previous level ($\ell = 0$ is the initialization step). Let
  $P$ be $A$'s parent and $C$ be a child of $B$. Throughout, we shall
  use the notation $C\succ B$ when $C$ is a child of $B$. At level
  $\ell-1$, the expansion coefficients $\{\delta^{P C}_{s}\}_{1\le
    s\le r_\eps}$ of $\{u^{C}(x)\}_{x\in P}$ are readily available
  and we have
  \[
  \left| u^{C}(x) - \sum_{s=1}^{r_\eps} \alpha^{PC}_{s}(x) \delta^{PC}_{s} \right| \le \left( \sum_{\xi\in C} |g(\xi)| \right) \eps,
  \quad \forall x\in P.
  \]
  Since $u^B(x) = \sum_{C\succ B} u^{C}(x)$, the previous inequality
  implies that
  \[
  \left| u^B(x) - \sum_{C\succ B} \sum_{s=1}^{r_\eps} \alpha^{PC}_{s}(x) \delta^{PC}_{s} \right| \le \left( \sum_{\xi\in B} |g(\xi)| \right) \eps,
  \quad \forall x\in P.
  \]
  Since $A \subset P$, the above approximation is of course true for
  any $x \in A$. However, since $\ell_A + \ell_B = L$, the sequence of
  restricted potentials $\{u^B(x)\}_{x\in A}$ also has a low-rank
  approximation of size $r_\eps$, namely,
  \[
  \left| u^B(x) - \sum_{t=1}^{r_\eps} \alpha^{AB}_t(x) \delta^{AB}_t \right| \le \left( \sum_{\xi\in B} |g(\xi)| \right) \eps,
  \quad \forall x\in A.
  \]
  Combining the last two approximations, we obtain that
  $\{\delta^{AB}_t\}_{1\le t\le r_\eps}$ should obey
  \begin{equation}
    \sum_{t=1}^{r_\eps} \alpha^{AB}_t(x) \delta^{AB}_t \approx
    \sum_{C\succ B} \sum_{s=1}^{r_\eps} \alpha^{PC}_{s}(x) \delta^{PC}_{s}, \quad \forall x\in A.
    \label{eqn:bf2}
  \end{equation}
  This is an over-determined linear system for
  $\{\delta^{AB}_t\}_{1\le t\le r_\eps}$ when
  $\{\delta^{PC}_{s}\}_{1\le s\le r_\eps,C\succ B}$ are available.
  Instead of computing $\{\delta^{AB}_t\}_{1\le t\le r_\eps}$ with a
  least-square method, the butterfly algorithm typically uses an
  efficient linear transformation approximately mapping
  $\{\delta^{PC}_{s}\}_{1\le s\le r_\eps,C\succ B}$ into
  $\{\delta^{AB}_t\}_{1\le t\le r_\eps}$. The actual implementation of
  this step is very much application-dependent.

\item Finally, $\ell = L$ and set $B$ to be the root node of
  $T_\Omega$. For each leaf box $A \in T_X$, use the constructed
  expansion coefficients $\{\delta^{AB}_t\}_{1\le t\le r_\eps}$ to
  evaluate $u(x)$ for each $x \in A$,
  \begin{equation}
    u(x) = \sum_{t=1}^{r_\eps} \alpha^{AB}_t (x) \delta^{AB}_t.
    \label{eqn:bf3}
  \end{equation}
\end{enumerate}

A schematic illustration of this algorithm is provided in Figure
\ref{fig:domain-tree-BF}. We would like to emphasize that the
strict balance between the levels of the target boxes $A$ and source
boxes $B$ maintained throughout this procedure is the key to obtain
the accurate low-rank separated approximations.

\input{fig-domain-tree-BF}

\section{Low-rank approximations}
\label{sec:LRapproximations}

In this section, the set $X$ and $\Omega$ refer to the sets defined in
\eqref{eqn:X} and \eqref{eqn:Omega}. In order to apply the algorithm
in Section \ref{eqn:FIOdiscrete}, one would require the existence
of the following low-rank separated representation
\[
e^{2\pi \i \Phi(x,\xi)} \approx \sum_{t=1}^{r_\eps} \alpha^{AB}_t(x) \beta^{AB}_t(\xi)
\]
for any pair of boxes $A$ and $B$ such that $\ell_A + \ell_B =
L$. However, this is not true for a general FIO kernel
$e^{2\pi\i\Phi(x,\xi)}$ due to the singularity of $\Phi(x,\xi)$ at the
origin $\xi=0$, i.e., when the square $B$ in $\Omega$ is close to the
origin of the frequency domain. However, if the frequency domain $B$
is well-separated from the origin $\xi=0$ in a relative sense, one can
prove a low-rank separated representation.

  In order to make it more precise, for two given squares $A\subset X$
  and $B\subset \Omega$, we introduce a new function called the
  residue phase function
\begin{equation}
  R^{AB}(x,\xi) := \Phi(x,\xi)-\Phi(c_A,\xi)-\Phi(x,c_B) +
  \Phi(c_A,c_B),
  \label{eqn:RAB}
\end{equation}
where $c_A$ and $c_B$ are the centers of $A$ and $B$
respectively. Using this new definition, the kernel can be written as
\begin{equation}
  e^{2\pi\i \Phi(x,\xi)} =
  e^{2\pi\i \Phi(c_A,\xi)} e^{2\pi\i \Phi(x,c_B)}
  e^{-2\pi\i\Phi(c_A,c_B)} e^{2\pi\i R^{AB}(x,\xi)}.
  \label{eqn:terms}
\end{equation}

\newcommand{\talpha}{\tilde{\alpha}}
\newcommand{\tbeta}{\tilde{\beta}}

\begin{theorem}
  \label{thm:hba}
  Suppose $\Phi(x,\xi)$ is a phase function that is real analytic for
  $x$ and $\xi$ away from $\xi=0$. There exists positive constants
  $\eps_0$ and $N_0$ such that the following is true.
  Let $A$ and $B$ be two
  squares in $X$ and $\Omega$, respectively, obeying $w_Aw_B\leq 1$
  and $\dist(B,0)\geq \frac{N}{4}$. For any positive $\eps\leq
  \eps_0$ and $N\geq N_0$, there exists an approximation
  \begin{equation*}
    \left|e^{2\pi \i R^{AB}(x,\xi)} -\sum^{r_\eps}_{t=1}
    \talpha^{AB}_t(x)\tbeta^{AB}_t(\xi)\right|
    \leq \eps
  \end{equation*}
  for $x\in A$ and $\xi\in B$ with $r_\eps\lesssim
  \log^4(\frac{1}{\eps})$. Moreover,
  \begin{itemize}
  \item when $w_B \le \sqrt{N}$, the functions
    $\{\tbeta^{AB}_t(\xi)\}_{1\le t\le r_\eps}$ can all be chosen as
    monomials in $(\xi-c_B)$ with a degree not exceeding a constant
    times $\log^2 (1/\eps)$,
  \item and when $w_A \le 1/\sqrt{N}$, the functions
    $\{\talpha^{AB}_t(x)\}_{1\le t\le r_\eps}$ can all be chosen as monomials in
    $(x-c_A)$ with a degree not exceeding a constant times $\log^2
    (1/\eps)$.
  \end{itemize}
\end{theorem}
In Theorem \ref{thm:hba}, $w_A$ and $w_B$ denote the side lengths of
$A$ and $B$, respectively; $\dist(B,0)$ denotes the distance between the
square $B$ and the origin $0$ in the frequency domain.  The distance is
given by $\dist(B,0) = \min_{\xi\in B}\norm{\xi-0}$.  Throughout this
paper, when we write $\O(\cdot)$, $\lesssim$ and $\gtrsim$, the
implicit constant is independent of $N$ and $\eps$.

\begin{proof}

Since $w_A w_B \le 1$, we either have $w_A \le 1/\sqrt{N}$ or $w_B
\le \sqrt{N}$ or both.

Let us first consider the case $w_B \le \sqrt{N}$. Then
\begin{align*}
  R^{AB}(x,\xi) & = \Phi(x,\xi) - \Phi(c_A,\xi) - \Phi(x,c_B)
      + \Phi(c_A,c_B)\\
  &  = \left[\Phi(x,\xi)-\Phi(c_A,\xi)\right]
      - \left[\Phi(x,c_B)-\Phi(c_A,c_B)\right]\\
  & = H(x,\xi) - H(x,c_B),
\end{align*}
where $H(x,\xi) := \Phi(x,\xi)-\Phi(c_A,\xi)$.
The function $R^{AB}(x,\xi)$ inherits the smoothness from
$\Phi(x,\xi)$. Applying the multi-variable Taylor expansion of degree
$k$ in $\xi$ centered at $c_B$ gives
\begin{equation}
  R^{AB}(x,\xi) =
  \sum_{1\le|i|<k} \frac{ \p^i_\xi H(x,c_B) }{i!} (\xi-c_B)^i
  + \sum_{|i|=k}   \frac{ \p^i_\xi H(x,\xi^*) }{i!} (\xi-c_B)^i,
  \label{eqn:RABexp}
\end{equation}
where $\xi^*$ is a point in the segment between $c_B$ and $\xi$.  Here
$i = (i_1,i_2)$ is a multi-index with $i! = i_1! i_2!$, and $|i| =
i_1+i_2$. Let us first choose the degree $k$ so that the second sum in
\eqref{eqn:RABexp} is bounded by $\eps/(4\pi)$. For each $i$ with $|i| = k$,
the definition of $H(x,\xi)$ gives
\[
\p^i_\xi H(x,\xi^*) = \sum_{|j|=1}
\p^j_x \p^i_\xi \Phi(x^*,\xi^*)(x-c_A)^j,
\]
for some point $x^*$ in the segment between $c_A$ and $x$. Using the
fact that $\Phi(x,\xi)$ is real-analytic over $|\xi|=1$ gives that
there exists a radius $R$ such that
\[
| \p_x^j \p_\xi^i \Phi(x,\xi) | \le C i! j!  \frac{1}{R^{|i+j|}}
= C i! j!  \frac{1}{R^{k+1}},
\]
for $\xi$ with $|\xi|=1$. Here the constant $C$ is independent of $k$.
Since $\Phi(x,\xi)$ is homogeneous of degree 1 in $\xi$, a scaling
argument shows that
\[
| \p_x^j \p_\xi^i \Phi(x^*,\xi^*) | \le C i! j!
\frac{1}{R^{k+1}  |\xi^*|^{k-1}}.
\]
Since $\dist(B,0)\ge N/4$ and $w_A w_B \le 1$, we have
\[
\left| \frac{\p^i_\xi H(x,\xi^*)}{i!} (\xi-c_B)^i \right| \le
\frac{ 2 C i! j!}{i!} \frac{1}{R^{k+1}  |\xi^*|^{k-1}} w_A  w_B^k
\le \frac{2C}{R^{k+1}} \left( \frac{4}{\sqrt{N}} \right)^{k-1}.
\]
Combining this with \eqref{eqn:RABexp} gives
\[
\left| R^{AB}(x,\xi) - \sum_{1\le|i|<k} \frac{ \p^i_\xi H(x,c_B) }{i!} (\xi-c_B)^i \right|
= \left| \sum_{|i|=k}   \frac{ \p^i_\xi H(x,\xi^*) }{i!} (\xi-c_B)^i \right|
\le \frac{2C(k+1)}{R^{k+1}} \left( \frac{4}{\sqrt{N}} \right)^{k-1}.
\]
Therefore, for a sufficient large $N_0(R)$, if $N>N_0(R)$, choosing $k
=k_\eps = O(\log(1/\eps))$ ensures that the difference is bounded by
$\eps/(4\pi)$.

The special case $k=1$ results in the following bound for
$R^{AB}(x,\xi)$
\[
|R^{AB}(x,\xi)| \le \frac{4C}{R^2}.
\]
To simplify the notation, we define
\[
R^{AB}_\eps(x,\xi) := \sum_{1\le|i|<k_\eps} \frac{ \p^i_\xi H(x,c_B)
}{i!}(\xi-c_B)^i,
\]
i.e., the first sum on the right hand side of \eqref{eqn:RABexp} with
$k=k_\eps$. The choice of $k_\eps$ together with \eqref{eqn:RABexp}
implies the bound
\[
|R^{AB}_\eps(x,\xi)| \le \frac{4C}{R^2} +\eps.
\]
Since $R_\eps^{AB}(x,\xi)$ is bounded, a direct application of Lemma
3.2 of \cite{FIO09} gives
\begin{equation}
  \left| e^{2\pi\i R_\eps^{AB}(x,\xi)} - \sum_{p=0}^{d_\eps}
  \frac{(2\pi\i R_\eps^{AB}(x,\xi))^p}{p!} \right| \le \eps/2,
  \label{eqn:expPAB}
\end{equation}
where $d_\eps = \O(\log(1/\eps))$. Since $R_\eps^{AB}(x,\xi)$ is a
polynomial in $(\xi-c_B)$, the sum in \eqref{eqn:expPAB} is also a
polynomial in $(\xi-c_B)$ with degree bounded by $k_\eps d_\eps =
\O(\log^2(1/\eps))$. Since our problem is in 2D, there are at most
$\O(\log^4(1/\eps))$ possible monomial in $(\xi-c_B)$ with degree
bounded by $k_\eps d_\eps$. Grouping the terms with the same
multi-index in $\xi$ results in an $\O(\log^4(1/\eps))$ term
$\eps$-accurate separated approximation for $e^{2\pi\i R_\eps^{AB}(x,\xi)}$ 
with the factors $\{\tbeta^{AB}_t(\xi)\}_{1\le t\le r_\eps}$ 
being monomials of $(\xi-c_B)$.

Finally, from the inequality $|e^{\imath a}-e^{\imath b}| \le |a-b|$,
it is clear
that a separated approximation for $e^{2\pi\i R^{AB}_\eps(x,\xi)}$
with accuracy $\eps/2$ is also one for $e^{2\pi\i R^{AB}(x,\xi)}$
with accuracy $\eps/2+\eps/2= \eps$. This completes the proof for the
case $w_B \le \sqrt{N}$.

The proof for the case $w_A \le 1/\sqrt{N}$ is similar. The only
difference is that we now group with
\[
R^{AB}(x,\xi) =
\left[\Phi(x,\xi)-\Phi(x,c_B)\right]-
\left[\Phi(c_A,\xi)-\Phi(c_A,c_B)\right]
\]
and apply the multivariable Taylor expansion in $x$ centered at $c_A$
instead. This results an $\O(\log^4(1/\eps))$ term $\eps$-accurate
separated approximation for $e^{2\pi\i R^{AB}(x,\xi)}$ with the
factors $\{\talpha^{AB}_t(x)\}_{1\le t\le r_\eps}$ being monomials of $(x-c_A)$.
\end{proof}

Though the above proof is constructive, it is cumbersome to construct
the separated approximation this way. On the other hand, the proof
shows that when $w_B \le \sqrt{N}$, the $\xi$-dependent factors in the
low-rank approximation of $e^{2\pi\i R^{AB}(x,\xi)}$ are all monomials
in $(\xi-c_B)$. Similarly, when $w_A \le 1/\sqrt{N}$, the
$x$-dependent factors are monomials in $(x-c_A)$.  This suggests to
use Chebyshev interpolation in $x$ when $w_A \le 1/\sqrt{N}$ and in
$\xi$ when $w_B \le \sqrt{N}$. For this purpose, we associate with
each box a Chebyshev grid as follows.

For a fixed integer $q$, the Chebyshev grid of order $q$ on
$[-1/2,1/2]$ is defined by
\[
\left\{ z_i = \frac{1}{2} \cos \left( \frac{i\pi}{q-1} \right) \right\}_{0\le i\le q-1}.
\]
A tensor-product grid {\em adapted to a square} with center $c$ and
side length $w$ is then defined via shifting and scaling as
\[
\{c + w (z_i, z_j)\}_{i,j=0,1,\ldots, q-1}
\]
In what follows, $M_t^B$ is the 2D Lagrange interpolation polynomial
on the Chebyshev grid adapted to the square $B$ (i.e., using $c=c_B$ and
$w=w_B$).

\begin{theorem}
  \label{thm:intp}
  Let $A$ and $B$ be as in Theorem \ref{thm:hba}. Then for any $\eps
  \le \eps_0$ and $ N \ge N_0$ where $\eps_0$ and $N_0$ are the
  constants in Theorem \ref{thm:hba}, there exists $q_\eps \lesssim
  \log^2 (1/\eps)$ such that
  \begin{itemize}
  \item when $w_B \le \sqrt{N}$, the Lagrange interpolation of
    $e^{2\pi\i R^{AB}(x,\xi)}$ in $\xi$ on a $q_\eps \times q_\eps$
    Chebyshev grid $\{g^B_t\}_{1\le t \le r_\eps}$ adapted to $B$
    obeys
    \begin{equation}
      \left| e^{2\pi\i R^{AB}(x,\xi)} - \sum_{t=1}^{r_\eps} e^{2\pi\i R^{AB}(x,g^B_t)} M^B_t(\xi) \right|
      \le \eps, \quad \forall x\in A, \forall \xi\in B,
      \label{eqn:intp1}
    \end{equation}
  \item when $w_A \le 1/\sqrt{N}$, the Lagrange interpolation of
    $e^{2\pi\i R^{AB}(x,\xi)}$ in $x$ on a $q_\eps \times q_\eps$
    Chebyshev grid $\{g^A_t\}_{1\le t\le r_\eps}$ adapted to $A$ obeys
    \begin{equation}
      \left| e^{2\pi\i R^{AB}(x,\xi)} - \sum_{t=1}^{r_\eps} M^A_t(x) e^{2\pi\i R^{AB}(g^A_t,\xi)} \right|
      \le \eps, \quad \forall x\in A, \forall \xi\in B.
      \label{eqn:intp2}
    \end{equation}
  \end{itemize}
  Both \eqref{eqn:intp1} and \eqref{eqn:intp2} provide a low-rank
  approximation with $r_\eps = q_\eps^2 \lesssim \log^4(1/\eps)$
  terms.
\end{theorem}
The proof for this follows exactly the one of Theorem 3.3 of
\cite{FIO09}.

Finally, we are ready to construct the low-rank approximation for the
kernel $e^{2\pi\i \Phi(x,\xi)}$, i.e.,
\begin{equation}
  e^{2\pi\i \Phi(x,\xi)} \approx  \sum^{r_\eps}_{t=1}\alpha^{AB}_t(x)\beta^{AB}_t(\xi).
  \label{eqn:PhiAB}
\end{equation}
When $w_B \le \sqrt{N}$, one multiply \eqref{eqn:intp1} with
$e^{2\pi\i \Phi(c_A,\xi)} e^{2\pi\i \Phi(x,c_B)}
e^{-2\pi\i\Phi(c_A,c_B)}$, which gives that $\forall x\in A, \forall
\xi\in B$
\[
\left|
e^{2\pi\i\Phi(x,\xi)} - \sum_{t=1}^{r_\eps} e^{2\pi\i\Phi(x,g^B_t)}
\left(
e^{-2\pi\i\Phi(c_A,g^B_t)} M^B_t(\xi) e^{2\pi\i\Phi(c_A,\xi)}
\right)
\right| \le \eps.
\]
In terms of the notations in \eqref{eqn:PhiAB}, the expansion functions
are given by
\begin{equation}
  \alpha^{AB}_t(x) = e^{2\pi\i\Phi(x,g^B_t)},\quad
  \beta^{AB}_t(\xi) = e^{-2\pi\i\Phi(c_A,g^B_t)} M^B_t(\xi) e^{2\pi\i\Phi(c_A,\xi)}, \quad 
  1\le t \le r_\eps.
\label{eqn:ab1}
\end{equation}
This is a special interpolant of the function $e^{2\pi\i\Phi(x,\xi)}$
in the $\xi$ variable, which pre-factors the oscillation, performs the
interpolation, and then remodulates the outcome.  When $w_A \le
1/\sqrt{N}$, multiply \eqref{eqn:intp2} with $e^{2\pi\i \Phi(c_A,\xi)}
e^{2\pi\i\Phi(x,c_B)} e^{-2\pi\i\Phi(c_A,c_B)}$ and obtain that
$\forall x\in A, \forall \xi\in B$
\[
\left|
  e^{2\pi\i\Phi(x,\xi)} - \sum_{t=1}^{r_\eps}
  \left(
  e^{2\pi\i\Phi(x,c_B)} M^A_t(x) e^{-2\pi\i\Phi(g^A_t,c_B)}
  \right)
  e^{2\pi\i\Phi(g^A_t,\xi)}
\right|
\le \eps.
\]
The expansion functions are now
\begin{equation}
  \alpha^{AB}_t(x) = e^{2\pi\i\Phi(x,c_B)} M^A_t(x) e^{-2\pi\i\Phi(g^A_t,c_B)},\quad
  \beta^{AB}_t(\xi) = e^{2\pi\i\Phi(g^A_t,\xi)}, \quad
  1\le t \le r_\eps.
  \label{eqn:ab2}
\end{equation}

Due to the presence of the demodulation and remodulation steps in the
definitions \eqref{eqn:ab1} and \eqref{eqn:ab2}, we refer to them as
{\em oscillatory Chebyshev interpolations}.

\section{Multiscale Butterfly Algorithm}
\label{sec:HBA}

In this section, we combine the low-rank approximations described in
Section \ref{sec:LRapproximations} with the butterfly algorithm in
Section \ref{sec:butterfly}. Due to the restriction on the distance
between $B$ and the origin, we decompose \eqref{eqn:FIOdiscrete} into
a multiscale summation
\begin{equation}
  \label{eqn:HFIO2}
  (L f)(x)=  \sum_{\xi \in \Omega_d}  e^{2\pi \i \Phi(x,\xi)}\widehat{f}(\xi)
  + \sum_j \sum_{\xi \in \Omega_j}  e^{2\pi \i \Phi(x,\xi)}\widehat{f}(\xi),
\end{equation}
where
\[
\Omega_j = \left\{(n_1,n_2):  \frac{N}{2^{j+1}} < \max(|n_1|,|n_2|) \le \frac{N}{2^j} \right\} \cap \Omega
\]
for $j=1,\dots,\log N - s$, $s$ is a constant, and
$\Omega_d=\Omega\setminus \cup_j\Omega_j$.

The term of $\Omega_d$ can be evaluated directly since $|\Omega_d| =
\O(1)$. Let us now fix an $\Omega_j$. Since any square $B$ in
$\Omega_j$ always stays away from the origin, the results in Section
\ref{sec:LRapproximations} applies to the term for $\Omega_j$ in
\eqref{eqn:HFIO2}. Therefore, the butterfly algorithm as described in
Section \ref{sec:butterfly} can be adapted to evaluate
\begin{equation*}
  \sum_{\xi \in \Omega_j}
  e^{2\pi \i \Phi(x,\xi)}\widehat{f}(\xi)
\end{equation*}
for the Cartesian domains $X$ and $\Omega_j$. In contrast to the
polar Butterfly algorithm that works in the polar coordinates for
$\Omega$, we refer to this one as {\em the Cartesian
  butterfly algorithm}.

\subsection{Cartesian butterfly algorithm}

To make it more explicit, let us first consider the interaction
between $(X,\Omega_1)$, with the low-rank approximation implemented
using the oscillatory Chebyshev interpolation discussed in Section
\ref{sec:LRapproximations}.
\begin{enumerate}
\item{\emph{Preliminaries.}} Construct two quadtrees $T_X$ and
  $T_{\Omega_1}$ for $X$ and $\Omega_1$ by uniform hierarchical
  partitioning. Let $b$ be a constant greater than or equal to $4$ and
  define $N_1 = N$.
\item{\emph{Initialization.}} For each square $A\in T_X$ of width $1/b$
  and each square $B\in T_{\Omega_1}$ of width $b$, the low-rank
  approximation functions are
  \begin{eqnarray}
    \label{eqn:intxi}
    \alpha_t^{AB}(x)=e^{2\pi \i \Phi(x,g_t^B)},\quad
    \beta_t^{AB}(\xi)=e^{-2\pi\i\Phi(c_A,g_t^B)}M_t^B(\xi)e^{2\pi\i\Phi(c_A,\xi)},\quad
    1\le t\le r_\eps.
  \end{eqnarray}
  Hence, we can define the expansion weights $\{\delta_t^{AB}\}_{1\le
    t\le r_\eps}$ with
  \begin{equation}
    \label{eqn:1}
    \delta_t^{AB} := \sum_{\xi\in B} \beta_t^{AB}(\xi)\widehat{f}(\xi)=
    e^{-2\pi\i \Phi(c_A,g_t^B)} \sum_{\xi\in B}\left( M_t^B(\xi)
    e^{2\pi\i\Phi(c_A,\xi)}\widehat{f}(\xi)\right).
  \end{equation}
\item{\emph{Recursion.}} Go up in tree $T_{\Omega_1}$ and down in tree
  $T_X$ at the same time until we reach the level such that $w_B =
  \sqrt{N_1}$. At each level, visit all the pairs $(A,B)$.  We apply
  the Chebyshev interpolation in variable $\xi$ and still define the
  approximation functions given in \eqref{eqn:intxi}.  Let
  $\{\delta_{s}^{P C}\}_{1\le s\le r_\eps}$ denote the expansion
  coefficients available in previous steps, where $P$ is $A$'s parent,
  $C$ is a child of $B$, and $s$ indicates the Chebyshev grid points
  in previous domain pairs. We define the new expansion coefficients
  $\{\delta_t^{AB}\}_{1\le t\le r_\eps}$ as
  \begin{equation}
    \label{eqn:2}
    \delta_t^{AB} := e^{-2\pi\i\Phi(c_A,g_t^B)}\sum_{C\succ B}\sum_{s=1}^{r_\eps} M_t^B(g_{s}^{C})
    e^{2\pi\i\Phi(c_A,g_{s}^{C})}\delta_{s}^{PC},
  \end{equation}
  where we recall that the notation $C\succ B$ means that $C$ is a
  child of $B$.
\item{\emph{Switch.}} For the levels visited, the Chebyshev
  interpolation is applied in variable $\xi$, while the interpolation
  is applied in variable $x$ for levels $l>\log(N_1)/2$.  Hence, we
  are switching the interpolation method at this step.  Now we are
  still working on level $l=\log(N_1)/2$ and the same domain pairs
  $(A,B)$ in the last step.  Let $\delta_s^{AB}$ denote the expansion
  weights obtained by Chebyshev interpolation in variable $\xi$ in the
  last step.  Correspondingly, $\{g_s^B\}_s$ are the grid points in
  $B$ in the last step.  We take advantage of the interpolation in
  variable $x$ in $A$ and generate grid points $\{g_t^A\}_{1\le t\le
    r_\eps}$ in $A$.  Then we can define new expansion weights
  \[
  \delta_t^{AB} := \sum_{s=1}^{r_\eps} e^{2\pi\i\Phi(g_t^A,g_s^B)}\delta_s^{AB}.
  \]
\item{\emph{Recursion.}} Go up in tree $T_{\Omega_1}$ and down in tree
  $T_X$ at the same time until we reach the level such that
  $w_B=N_1/b$. We construct the approximation functions by Chebyshev
  interpolation in variable $x$ as follows:
  \begin{eqnarray}
    \label{eqn:intx}
    \alpha_t^{AB}(x)= e^{2\pi\i\Phi(x,c_B)}M_t^A(x) e^{-2\pi\i\Phi(g_t^A,c_B)},
    &\beta_t^{AB}(\xi)=e^{2\pi\i\Phi(g_t^A,\xi)}.
  \end{eqnarray}
  We define the new expansion coefficients $\{\delta_t^{AB}\}_{1\le
    t\le r_\eps}$ as
  \begin{equation}
    \label{eqn:3}
    \delta_t^{AB} := \sum_{C\succ B} e^{2\pi\i\Phi(g_t^A,c_{C})}
    \sum_{s=1}^{r_\eps}\left(M_{s}^{P}(g_{t}^{A})
    e^{-2\pi\i\Phi(g_{s}^{P},c_{C})}\delta_{s}^{PC}\right),
  \end{equation}
  where again $P$ is $A$'s parent and $C$ is a child box of $B$.
\item{\emph{Termination.}} Finally, we reach the level that
  $w_B=N_1/b$. For each $B$ on this level and for each square $A\in
  T_X$ of width $b/N_1$, we apply the approximation functions given by
  \eqref{eqn:intx} and obtain
  \begin{equation}
  \label{eqn:sub}
    u^{B}(x) :=
    e^{2\pi\i\Phi(x,c_B)} \sum_{t=1}^{r_\eps}\left(M_t^A(x)e^{-2\pi\i\Phi(g_t^A,c_B)} \delta_t^{AB} \right)
  \end{equation}
  for each $x\in A$. Finally, summing over all $B$ on this level, we have
   \begin{equation}
    \label{eqn:4}
  u^{\Omega_1}(x) :=\sum_{B} u^{B}(x)
  \end{equation}
   for each $x\in A$. 
\end{enumerate}
We would like to emphasize that the center part of the tree
$T_{\Omega_j}$ is always empty since $\Omega_j$ is a
corona. Accordingly, the algorithm skips this empty part.

For a general $\Omega_j$, the interaction between $(X,\Omega_j)$
follows a similar algorithm, except that we replace $\Omega_1$ with
$\Omega_j$, $u^{\Omega_1}(x)$ with $u^{\Omega_j}(x)$,
$N_1$ with $N_j= N/2^{j-1}$, and stop at the level that
$w_B = N_j/b$.

Finally, \eqref{eqn:HFIO2} is evaluated via
\begin{equation}
(Lf)(x) = u^{\Omega_d}(x) + \sum_j u^{\Omega_j}(x).
\end{equation}

\subsection{Complexity analysis}

The cost of evaluating the term of $\Omega_d$ takes at most $\O(N^2)$
steps since $|\Omega_d|=\O(1)$. Let us now consider the cost of the
terms associated with $\{\Omega_j\}$.

For the interaction between $X$ and $\Omega_1$, the computation
consists of two parts: the recursive evaluation of $\{\delta^{AB}_t\}$
and the final evaluation of $u^{\Omega_1}(x)$. The recursive part takes
$\O(q^3 N^2 \log N)$ since there are at most $\O(N^2 \log N)$ pairs of
squares $(A,B)$ and the evaluation of $\{\delta^{AB}_t\}$ for each
pair takes $\O(q^3)$ steps via dimension-wise Chebyshev
interpolation. The final evaluation of $u^{\Omega_1}(x)$ clearly takes
$\O(q^2 N^2)$ steps as we spend $\O(q^2)$ on each point $x\in X$.

For the interaction between $X$ and $\Omega_j$, the analysis is
similar. The recursive part takes now $\O(q^3 N_j^2 \log N_j)$ steps
(with $N_j=N/2^{j-1}$) as there are at most $\O(N_j^2\log N_j)$ pairs
of squares involved. The final evaluation still takes $\O(q^2 N^2)$
steps.

Summing these contributions together results in the total
computational complexity
\[
\O(q^3 N^2 \log N) + \O(q^2 N^2\log N) = \O(q^3 N^2 \log N) =
\O(r_\eps^{3/2} N^2 \log N).
\]
The multiscale butterfly algorithm is also highly efficient in terms
of memory as the Cartesian butterfly algorithm is applied sequentially
to evaluate \eqref{eqn:sub} for each $\Omega_j$. The overall memory
complexity is $O(\frac{N^2}{b^2})$, only $\frac{1}{b^2}$ of that the
original Cartesian butterfly algorithm.

\section{Numerical results}
\label{sec:results}
This section presents several numerical examples to demonstrate the
effectiveness of the multiscale butterfly algorithm introduced above.
In truth, FIOs usually have non-constant amplitude functions.
Nevertheless, the main computational difficulty is the oscillatory
phase term.  We refer to \cite{FIO09} for detailed fast algorithms to
deal with non-constant amplitude functions.  Our MATLAB implementation
can be found on the authors' personal homepages.  The numerical
results were obtained on a desktop with a 3.5 {GHz CPU} and 32 {GB} of
memory.  Let $\{u^d(x),x\in X\}$, $\{u^m(x),x\in X\}$ and
$\{u^p(x),x\in X\}$ be the results of a discrete FIO computed by a
direct matrix-vector multiplication, the multiscale butterfly
algorithm and the polar butterfly algorithm~\cite{FIO09},
respectively.  To report on the accuracy, we randomly select a set $S$
of 256 points from $X$ and evaluate the relative errors of the
multiscale butterfly algorithm and the polar butterfly algorithm by
\begin{equation}
\eps^m=\sqrt{\frac{\sum_{x\in S}\lvert u^d(x)-u^m(x)\rvert^2}
{\sum_{x\in S}\lvert u^d(x)\rvert^2}}
\text{ and }
\eps^p=\sqrt{\frac{\sum_{x\in S}\lvert u^d(x)-u^p(x)\rvert^2}
{\sum_{x\in S}\lvert u^d(x)\rvert^2}}.
\end{equation}

According to the description of the multiscale butterfly algorithm in
Section \ref{sec:HBA}, we recursively divide $\Omega$ into $\Omega_j,
j=1,2,\dots,\log N-s$, where $s$ is $5$ in the following examples.
This means that the center square $\Omega_d$ is of size $2^5\times 2^5$
and the interaction from $\Omega_d$ is evaluated via a direct matrix-vector
multiplication.  Suppose $q_\eps$ is the number of Chebyshev points in each
dimension.  There is no sense to use butterfly algorithms to construct
$\{\delta_t^{AB}\}$ when the number of points in $B$ is fewer than
$q_\eps^2$. Hence, the recursion step in butterfly algorithms starts from
the squares $B$ that are a couple of levels away from the bottom of
$T_\Omega$ such that each square contains at least $q_\eps^2$ points.
Similarly, the recursion stops at the squares in $T_X$ that are the same
number of levels away from the bottom.  In the following examples, we
start from level $\log N-3$ and stop at level 3
(corresponding to $b=2^3$ defined in Section \ref{sec:HBA})
which matches with $q_\eps$ (4 to 11).

In order to make a fair comparison, we compare the MATLAB versions of
the polar butterfly algorithm and the multiscale butterfly algorithm.
Hence, the running time of the polar butterfly algorithm here is
slower than the one in \cite{FIO09}, which was implemented in C++.

\textbf{Example 1.} The first example is a generalized Radon transform
whose kernel is given by
\begin{equation}
\label{eqn:ex1}
\begin{split}
\Phi(x,\xi) &= x\cdot \xi + \sqrt{c_1^2(x)\xi_1^2+c_2^2(x)\xi_2^2},\\
c_1(x)&=(2+\sin(2\pi x_1)\sin(2\pi x_2))/3,\\
c_2(x)&=(2+\cos(2\pi x_1)\cos(2\pi x_2))/3.
\end{split}
\end{equation}
We assume the amplitude of this example is a constant $1$.
Now the FIO models an integration
over ellipses where $c_1(x)$ and $c_2(x)$ are the axis lengths
of the ellipse centered at the point $x\in X$.
Table~\ref{tab:bfio} summarize
the results of this example given by
the polar butterfly algorithm and the multiscale butterfly algorithm.


\begin{table}[htp]
\centering
\begin{tabular}{|rcc|rcc|c|}
\hline
\multicolumn{3}{|c|}{Multiscale Butterfly} & \multicolumn{3}{c|}{Polar Butterfly} & \\
\hline
     $N,q_\eps$ & $\eps^m$ &$T_m(sec)$ &
     $N,q_\eps$ & $\eps^p$ &$T_p(sec)$ &$T_p/T_m$\\
\hline
   256,5 & 7.89e-02 & 6.96e+01 &    256,5 & 4.21e-02 & 4.84e+02 & 6.96e+00 \\
   512,5 & 9.01e-02 & 3.62e+02 &    512,5 & 5.54e-02 & 2.34e+03 & 6.46e+00 \\
  1024,5 & 9.13e-02 & 1.81e+03 &   1024,5 & 4.26e-02 & 1.14e+04 & 6.31e+00 \\
  2048,5 & 9.47e-02 & 8.79e+03 &   2048,5 & - & - & -  \\
\hline
   256,7 & 6.95e-03 & 8.20e+01 &    256,7 & 5.66e-03 & 5.97e+02 & 7.28e+00 \\
   512,7 & 8.43e-03 & 4.16e+02 &    512,7 & 5.89e-03 & 2.82e+03 & 6.79e+00 \\
  1024,7 & 8.45e-03 & 2.03e+03 &   1024,7 & 4.84e-03 & 1.35e+04 & 6.64e+00 \\
  2048,7 & 8.42e-03 & 1.04e+04 &   2048,7 & - & - & -  \\
\hline
   256,9 & 3.90e-04 & 1.10e+02 &    256,9 & 8.25e-04 & 7.74e+02 & 7.04e+00 \\
   512,9 & 3.42e-04 & 5.39e+02 &    512,9 & 6.78e-04 & 3.57e+03 & 6.61e+00 \\
  1024,9 & 7.61e-04 & 2.74e+03 &   1024,9 & 4.18e-04 & 1.67e+04 & 6.09e+00 \\
  2048,9 & 4.82e-04 & 1.25e+04 &   2048,9 & - & - & -  \\
\hline
   256,11 & 2.15e-05 & 1.84e+02 &    256,11 & 3.69e-05 & 1.15e+03 & 6.27e+00 \\
   512,11 & 1.89e-05 & 8.60e+02 &    512,11 & 5.53e-05 & 5.10e+03 & 5.93e+00 \\
  1024,11 & 1.96e-05 & 4.27e+03 &    1024,11 & 2.042e-05 & 2.30e+04 &  5.39e+00 \\
  2048,11 & 1.50e-05 & 1.82e+04 &   2048,11 & - & - & - \\
\hline
\end{tabular}
\caption{Comparison of the multiscale butterfly algorithm
    and the polar butterfly algorithm
    for the phase function in \eqref{eqn:ex1}.
    $T_{m}$ is the running time of the multiscale butterfly algorithm;
    $T_{a}$ is the running time of the polar butterfly algorithm;
    and $T_m/T_p$ is the speedup factor.}
\label{tab:bfio}
\end{table}

\textbf{Example 2.} Next, we provide an FIO example with a smooth
amplitude function,
\begin{equation}
\label{eqn:ex2}
u(x)=\sum_{\xi\in\Omega}a(x,\xi)e^{2\pi\i\Phi(x,\xi)}\widehat{f}(\xi),
\end{equation}
where the amplitude and phase functions are given by
\begin{eqnarray*}
&a(x,\xi)=(J_0(2\pi \rho(x,\xi))+
\i Y_0(2\pi \rho(x,\xi)))e^{-\pi\i \rho(x,\xi)},\\
&\Phi(x,\xi)=x\cdot \xi + \rho(x,\xi),\\
&\rho(x,\xi)=\sqrt{c_1^2(x)\xi_1^2+c_2^2(x)\xi_2^2},\\
&c_1(x)=(2+\sin(2\pi x_1)\sin(2\pi x_2))/3,\\
&c_2(x)=(2+\cos(2\pi x_1)\cos(2\pi x_2))/3.
\end{eqnarray*}
Here, $J_0$ and $Y_0$ are Bessel functions of the first and second
kinds.  We refer to \cite{FIO07} for more details of the derivation of
these formulas.  As discussed in \cite{FIO09}, we compute the low rank
approximation of the amplitude functions $a(x,\xi)$ first:
\[
a(x,\xi)\approx \sum_{t=1}^{s_\eps}g_{t}(x)h_{t}(\xi).
\]
In the second step, we apply the multiscale butterfly algorithm to
compute
\begin{equation*}
u_{t}(x)=\sum_{\xi\in\Omega}e^{2\pi\i\Phi(x,\xi)}
\widehat{f}(\xi)h_{t}(\xi),
\end{equation*}
and sum up all $g_t(x)u_t(x)$ to evaluate
\begin{equation*}
u(x)=\sum_{t}g_t(x)u_t(x).
\end{equation*}
Table~\ref{tab:hbfio2} summarizes the results of this example given by
the direct method and the multiscale butterfly algorithm.

\begin{table}[htp]
\centering
\begin{tabular}{|rcccc|}
\hline
$N,q_\eps$ & $\eps^m$ &$T_d(sec)$&$T_m(sec)$&$T_d/T_m$ \\
\hline
   256,7 & 5.10e-03 & 3.78e+03 & 6.07e+02 & 6.23e+00 \\
   512,7 & 7.29e-03 & 3.71e+04 & 3.50e+03 & 1.06e+01 \\
  1024,7 & 6.16e-03 & 6.42e+05 & 1.70e+04 & 3.77e+01 \\
\hline
   256,9 & 4.49e-04 & 2.34e+03 & 7.88e+02 & 2.97e+00 \\
   512,9 & 4.04e-04 & 3.66e+04 & 4.64e+03 & 7.90e+00 \\
  1024,9 & 3.88e-04 & 6.21e+05 & 2.17e+04 & 2.86e+01 \\
\hline
   256,11 & 1.86e-05 & 2.48e+03 & 1.33e+03 & 1.86e+00 \\
   512,11 & 1.80e-05 & 3.60e+04 & 6.94e+03 & 5.18e+00 \\
  1024,11 & 2.39e-05 & 5.96e+05 & 2.83e+04 & 2.11e+01 \\
\hline
\end{tabular}
\caption{Numerical results given by the multiscale butterfly algorithm
  for the FIO in \eqref{eqn:ex2}.  $T_{d}$ is the running time of the
  direct evaluation; $T_{m}$ is the running time of the multiscale
  butterfly algorithm; and $T_d/T_m$ is the speedup factor.}
\label{tab:hbfio2}
\end{table}

Note that the accuracy of the multiscale butterfly algorithm is well
controlled by the number of Chebyshev points $q_\eps$. This indicates that our
algorithm is numerically stable.  Another observation is that the
relative error improves on average by a factor of 12 every time $q_\eps$ is
increased by a factor of 2. As we can see in those tables, for a fixed
kernel and a fixed $q_\eps$, the accuracy is almost independent of $N$.
Hence, in practical applications, one can increase the value of $q_\eps$
until a desired accuracy is reached in the problem with a small $N$.
In the comparison in Table \ref{tab:bfio}, the multiscale butterfly
algorithm and the polar butterfly algorithm use $q_\eps=\{5,7,9,11\}$ and
achieve comparable accuracy.  Meanwhile, as we observed from Table
\ref{tab:bfio}, the relative error decreasing rate of
the multiscale butterfly algorithm is larger than the decreasing rate
of the polar butterfly algorithm.
This means if a high accuracy is desired,
the multiscale butterfly algorithm requires a smaller $q_\eps$ to achieve
it comparing to the polar butterfly algorithm.

The second concern about the algorithm is the asymptotic complexity.
From the $T_m$ column of Table \ref{tab:bfio} and \ref{tab:hbfio2}, we
see that $T_m$ almost quadrupled when the problem size doubled under
the same $q_\eps$.  According to this, we are convinced that the empirical
running time of the multiscale butterfly algorithm follows the
$\O\left(N^2\log N\right)$ asymptotic complexity.  Note that the
speedup factor over the polar butterfly algorithm is about $6$ and the
multiscale butterfly algorithm obtains better accuracy.  This makes
the multiscale butterfly algorithm quite attractive to practitioners
who are interested in evaluating an FIO with a large $N$.

\textbf{Example 3.} Extending the multiscale butterfly algorithm to
higher dimensions is straightforward.  There are two main
modifications: higher dimensional multiscale domain decomposition and
Chebyshev interpolation.  In three dimensions, the frequency domain is
decomposed into cubic shells instead of coronas.  The
kernel interpolation is applied on a three dimensional Chebyshev
grids.  We apply our three-dimensional multiscale butterfly algorithm
to a simple example integrating over spheres with different radii.  We
assume a constant amplitude function and the kernel function is given
by
\begin{equation}\label{eqn:ex3}
\Phi(x,\xi) = x\cdot \xi +c(x)\sqrt{\xi_1^2+\xi_2^2},\quad
c(x) = (3+\sin(2\pi x_1)\sin(2\pi x_2)\sin(2\pi x_3))/4.
\end{equation}
Table \ref{tab:hbfio3} summarizes the results of this example given by
the direct method and the multiscale butterfly algorithm.

\begin{table}[htp]
\centering
\begin{tabular}{|rcccc|}
\hline
$N,q_\eps$ & $\eps^m$ &$T_d(sec)$&$T_m(sec)$&$T_d/T_m$ \\
\hline
   64,5 & 9.41e-02 & 1.82e+04 & 2.50e+03 & 7.31e+00 \\
   128,5 & 7.57e-02 & 6.21e+05 & 2.42e+04 & 2.57e+01 \\
   256,5 & 8.23e-02 & 3.91e+07 & 2.35e+05 & 1.66e+02 \\
\hline
    64,7 & 1.20e-02 & 1.83e+04 & 7.32e+03 & 2.50e+00 \\
   128,7 & 1.03e-02 & 6.03e+05 & 4.48e+04 & 1.35e+01 \\
   256,7 & 8.13e-03 & 4.39e+07 & 3.81e+05 & 1.15e+02 \\
\hline
\end{tabular}
\caption{Numerical results given by
    the multiscale butterfly algorithm for the phase function in \eqref{eqn:ex3}.}
\label{tab:hbfio3}
\end{table}

\section{Conclusion}
\label{sec:conclusion}
A simple and efficient multiscale butterfly algorithm for evaluating
FIOs is introduced in this paper.  This method hierarchically
decomposes the frequency domain into multiscale coronas in order to
avoid possible singularity of the phase function $\Phi(x,\xi)$ at
$\xi=0$.  A Cartesian butterfly algorithm is applied to evaluate the
FIO over each corona. Many drawbacks of the original butterfly
algorithm based on a polar-Cartesian transform in \cite{FIO09} can be
avoided. The new multiscale butterfly algorithm has an $\O(N^2\log N)$
operation complexity with a smaller pre-factor, while keeping the same
$\O(N^2)$ memory complexity.

{\bf Acknowledgments.}  This work was partially supported by the
National Science Foundation under award DMS-1328230 and the
U.S. Department of Energy’s Advanced Scientific Computing Research
program under award DE-FC02-13ER26134/DE-SC0009409.

\bibliographystyle{abbrv} \bibliography{ref}

\end{document}

%% file: fig-domain-decomp.tex
\begin{figure}[htb]
\centering
\begin{tikzpicture}[scale=1]
\coordinate (T) at (0,-2.7);

\fill[lightgray] (-2,-2) rectangle (2,2);
\draw[black] (-2,-2) rectangle (2,2);
\fill[white] (-1,-1) rectangle (1,1);
\draw[black] (-1,-1) rectangle (1,1);

\draw ($(T)$) node[rectangle] {$\Omega_1$};

\coordinate (A) at (4,0);
\fill[lightgray] ($(A)+(-1,-1)$) rectangle ($(A)+(1,1)$);
\draw[black] ($(A)+(-1,-1)$) rectangle ($(A)+(1,1)$);
\fill[white] ($(A)+(-.5,-.5)$) rectangle ($(A)+(.5,.5)$);
\draw[black] ($(A)+(-.5,-.5)$) rectangle ($(A)+(.5,.5)$);

\draw ($(A)+(T)$) node[rectangle] {$\Omega_2$};

\coordinate (A) at (7,0);
\draw ($(A)$) node[rectangle] {$\cdots$};
\draw ($(A)+(T)$) node[rectangle] {$\cdots$};

\coordinate (A) at (9,0);
\fill[lightgray] ($(A)+(-.25,-.25)$) rectangle ($(A)+(.25,.25)$);
\draw[black] ($(A)+(-.25,-.25)$) rectangle ($(A)+(.25,.25)$);
\fill[white] ($(A)+(-.125,-.125)$) rectangle ($(A)+(.125,.125)$);
\draw[black] ($(A)+(-.125,-.125)$) rectangle ($(A)+(.125,.125)$);

\draw ($(A)+(T)$) node[rectangle] {$\Omega_{\log N -s}$};

\coordinate (A) at (10.5,0);
\fill[lightgray] ($(A)+(-.125,-.125)$) rectangle ($(A)+(.125,.125)$);
\draw[black] ($(A)+(-.125,-.125)$) rectangle ($(A)+(.125,.125)$);
\draw ($(A)+(T)$) node[rectangle] {$\Omega_d$};

\end{tikzpicture}
\caption{This figure shows the frequency domain decomposition of $\Omega$.
    Each sub-domain $\Omega_j$, $j=1,\dots,\log N-s$,
    is a corona
    and $\Omega_d$ is a small square domain near the origin.}
\label{fig:domain-decomp}
\end{figure}
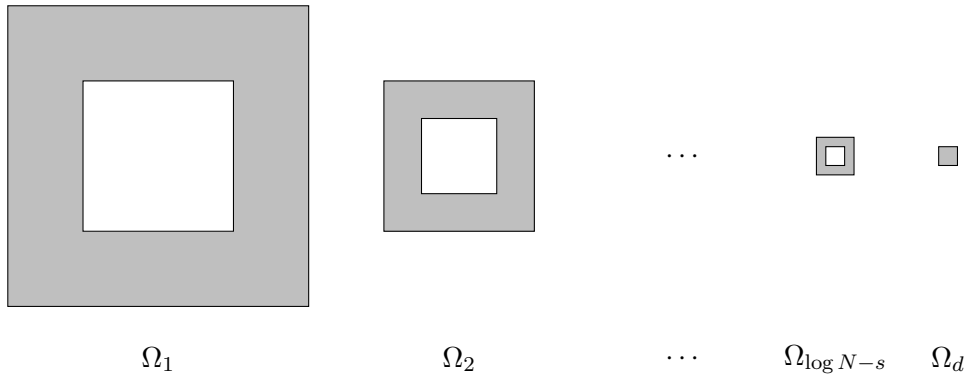

%% file: fig-domain-tree-BF.tex
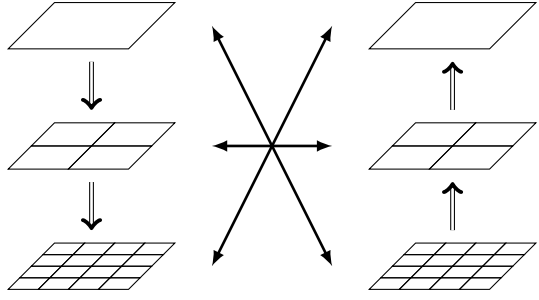
\begin{figure}
\centering
\begin{tikzpicture}[scale=0.4]

\draw [->,double,double distance=1pt] (2,-1.2,2) -> (2,-2.8,2);
\draw [->,double,double distance=1pt] (2,2.8,2) -> (2,1.2,2);

\draw (0,4,0) -- (0,4,4) -- (4,4,4) -- (4,4,0) -- cycle;

\draw (0,0,0) -- (0,0,2) -- (2,0,2) -- (2,0,0) -- cycle;
\draw (0,-4,0) -- (0,-4,1) -- (1,-4,1) -- (1,-4,0) -- cycle;
\draw (1,-4,1) -- (1,-4,2) -- (2,-4,2) -- (2,-4,1) -- cycle;
\draw (1,-4,0) -- (1,-4,1) -- (2,-4,1) -- (2,-4,0) -- cycle;
\draw (0,-4,1) -- (0,-4,2) -- (1,-4,2) -- (1,-4,1) -- cycle;

\draw (2,0,0) -- (2,0,2) -- (4,0,2) -- (4,0,0) -- cycle;
\draw (2,-4,0) -- (2,-4,1) -- (3,-4,1) -- (3,-4,0) -- cycle;
\draw (3,-4,1) -- (3,-4,2) -- (4,-4,2) -- (4,-4,1) -- cycle;
\draw (3,-4,0) -- (3,-4,1) -- (4,-4,1) -- (4,-4,0) -- cycle;
\draw (2,-4,1) -- (2,-4,2) -- (3,-4,2) -- (3,-4,1) -- cycle;

\draw (2,0,2) -- (2,0,4) -- (4,0,4) -- (4,0,2) -- cycle;
\draw (2,-4,2) -- (2,-4,3) -- (3,-4,3) -- (3,-4,2) -- cycle;
\draw (3,-4,3) -- (3,-4,4) -- (4,-4,4) -- (4,-4,3) -- cycle;
\draw (3,-4,2) -- (3,-4,3) -- (4,-4,3) -- (4,-4,2) -- cycle;
\draw (2,-4,3) -- (2,-4,4) -- (3,-4,4) -- (3,-4,3) -- cycle;

\draw (0,0,2) -- (0,0,4) -- (2,0,4) -- (2,0,2) -- cycle;
\draw (0,-4,2) -- (0,-4,3) -- (1,-4,3) -- (1,-4,2) -- cycle;
\draw (1,-4,3) -- (1,-4,4) -- (2,-4,4) -- (2,-4,3) -- cycle;
\draw (1,-4,2) -- (1,-4,3) -- (2,-4,3) -- (2,-4,2) -- cycle;
\draw (0,-4,3) -- (0,-4,4) -- (1,-4,4) -- (1,-4,3) -- cycle;

\draw [->,double,double distance=1pt] (2+12,-2.8,2) -> (2+12,-1.2,2);
\draw [->,double,double distance=1pt] (2+12,1.2,2) -> (2+12,2.8,2);

\draw (0+12,4,0) -- (0+12,4,4) -- (4+12,4,4) -- (4+12,4,0) -- cycle;

\draw (0+12,0,0) -- (0+12,0,2) -- (2+12,0,2) -- (2+12,0,0) -- cycle;
\draw (0+12,-4,0) -- (0+12,-4,1) -- (1+12,-4,1) -- (1+12,-4,0) -- cycle;
\draw (1+12,-4,1) -- (1+12,-4,2) -- (2+12,-4,2) -- (2+12,-4,1) -- cycle;
\draw (1+12,-4,0) -- (1+12,-4,1) -- (2+12,-4,1) -- (2+12,-4,0) -- cycle;
\draw (0+12,-4,1) -- (0+12,-4,2) -- (1+12,-4,2) -- (1+12,-4,1) -- cycle;

\draw (2+12,0,0) -- (2+12,0,2) -- (4+12,0,2) -- (4+12,0,0) -- cycle;
\draw (2+12,-4,0) -- (2+12,-4,1) -- (3+12,-4,1) -- (3+12,-4,0) -- cycle;
\draw (3+12,-4,1) -- (3+12,-4,2) -- (4+12,-4,2) -- (4+12,-4,1) -- cycle;
\draw (3+12,-4,0) -- (3+12,-4,1) -- (4+12,-4,1) -- (4+12,-4,0) -- cycle;
\draw (2+12,-4,1) -- (2+12,-4,2) -- (3+12,-4,2) -- (3+12,-4,1) -- cycle;

\draw (2+12,0,2) -- (2+12,0,4) -- (4+12,0,4) -- (4+12,0,2) -- cycle;
\draw (2+12,-4,2) -- (2+12,-4,3) -- (3+12,-4,3) -- (3+12,-4,2) -- cycle;
\draw (3+12,-4,3) -- (3+12,-4,4) -- (4+12,-4,4) -- (4+12,-4,3) -- cycle;
\draw (3+12,-4,2) -- (3+12,-4,3) -- (4+12,-4,3) -- (4+12,-4,2) -- cycle;
\draw (2+12,-4,3) -- (2+12,-4,4) -- (3+12,-4,4) -- (3+12,-4,3) -- cycle;

\draw (0+12,0,2) -- (0+12,0,4) -- (2+12,0,4) -- (2+12,0,2) -- cycle;
\draw (0+12,-4,2) -- (0+12,-4,3) -- (1+12,-4,3) -- (1+12,-4,2) -- cycle;
\draw (1+12,-4,3) -- (1+12,-4,4) -- (2+12,-4,4) -- (2+12,-4,3) -- cycle;
\draw (1+12,-4,2) -- (1+12,-4,3) -- (2+12,-4,3) -- (2+12,-4,2) -- cycle;
\draw (0+12,-4,3) -- (0+12,-4,4) -- (1+12,-4,4) -- (1+12,-4,3) -- cycle;

\draw [latex-latex,line width=1pt] (10,0,2) -> (6,0,2);
\draw [latex-latex,line width=1pt] (10,4,2) -> (6,-4,2);
\draw [latex-latex,line width=1pt] (10,-4,2) -> (6,4,2);

\end{tikzpicture}
\caption{Hierarchical domain trees of the 2D butterfly algorithm.
    Left: $T_X$ for the spatial domain $D_X$. 
    Right: $T_\Omega$ for the frequency domain $D_\Omega$. 
    The interactions between subdomains $A\subset D_X$
    and $B\subset D_\Omega$ are represented by left right arrow lines. }
\label{fig:domain-tree-BF}
\end{figure}